\documentstyle[11pt,amscd,amssymb]{amsart}
\theoremstyle{plain}
 \newtheorem{thm}{Theorem}[section]
 \newtheorem{lem}[thm]{Lemma}
 \newtheorem{prop}[thm]{Proposition}
 \newtheorem{cor}[thm]{Corollary}
 
\theoremstyle{definition}
 
\theoremstyle{remark}
 
 \newtheorem{ex}{Example}

\newcommand{ \Supp}{\operatorname{Supp}}
\newcommand{\Ext}{\operatorname{Ext}}
\newcommand{\Hom}{\operatorname{Hom}}

\newcommand{\rk}{\operatorname{rk}}

\newcommand{\NS}{\operatorname{NS}}
\newcommand{\coker}{\operatorname{coker}}
\newcommand{\Pic}{\operatorname{Pic}}
\newcommand{\ch}{\operatorname{ch}}
\newcommand{\td}{\operatorname{td}}

\newcommand{\Hilb}{\operatorname{Hilb}}

\newcommand{\IT}{\operatorname{IT}}
\newcommand{\WIT}{\operatorname{WIT}}

\font\b=cmr10 scaled \magstep5
\def\bigzerou{\smash{\lower1.7ex\hbox{\b 0}}}

\numberwithin{equation}{section}
\begin{document}

\title{
Some examples of isomorphisms induced by Fourier-Mukai functors
}
\author{K\={o}ta Yoshioka}
 \address{Max Planck Institut f\"{ur} Mathematik,
Gottfried Claren Str. 26, D-53225 Bonn, Germany $\&$
\newline
Department of mathematics, Faculty of Science, Kobe University,
Kobe, 657, Japan}
\email{yoshioka@@mpim-bonn.mpg.de}
 \maketitle

\section{Introduction}

In order to investigate sheaves on abelian varieties,
Mukai [Mu1] introduced a very powerful tool
called Fourier-Mukai functor.
As an application, Mukai ([Mu1], [Mu4]) computed some moduli spaces of stable sheaves on
abelian varieties.
Recently, Dekker [D] found some examples of isomorphisms of
moduli spaces of sheaves induced by Fourier-Mukai functors.
As an application, he proved that moduli spaces of sheaves 
on abelian surfaces are deformation equivalent to 
Hilbert schemes of points (see Theorem \ref{thm:another}).
 
Recently Fourier-Mukai functor was generalized to 
more general situations (e.g. [Br1], [Br2], [Mu6], [Mu7]).
Next task is to construct many examples of birational maps
of moduli spaces of sheaves.
In this note, 
we restrict ourselves to an abelian or a K3 surface of Picard
number 1 and 
give some examples 
of birational maps of moduli spaces of sheaves
induced by Fourier-Mukai functor (Theorem \ref{thm:fourier}).
More precisely, we shall find some examples of sheaves
which satisfy $\WIT_i$.
In general, Fourier-Mukai functor does not induce
isomorphisms of moduli spaces of sheaves.
Motivated by recent work of Markman [Mr],
we also consider the composition of Fourier-Mukai functor
and ``taking-dual'' functor.
Then we can get isomorphisms in these cases.

As an application, we get another proof of Dekker's result
(Theorem \ref{thm:another}).
Our condition and method are similar to [Y4].
In section 1, we shall treat original Fourier-Mukai functor.
In section 2, we shall explain how to generalize it 
to more general situations.

\vspace{1pc}

{\it Notation.}
Let $X$ be an abelian or a K3 surface defined over ${\Bbb C}$.
We set $H^{ev}(X,{\Bbb Z}):=\oplus_i H^{2i}(X,{\Bbb Z})$.
For $x \in H^{ev}(X,{\Bbb Z})$, $[x]_i$ denotes the $2i$-th
component of $x$.
Let $(H^{ev}(X,{\Bbb Z}),\langle \quad,\quad \rangle)$ be the Mukai lattice
of $X$.

Let ${\mathbf D}(X)$ be the derived category of $X$.
For $x \in {\mathbf D}(X)$,
\begin{equation}
\begin{split}
v(x): &=\ch(x)\sqrt{\td_X} \\
&=\ch(x)(1+\varepsilon \omega)
\in H^{ev}(X,{\Bbb Z})
\end{split}
\end{equation}
is the Mukai vector of $x$,
where $\varepsilon=0,1$ according as $X$ is an abelian surface
or a K3 surface,
and $\omega$ is the fundamental class of $X$.
Let $L$ be an ample line bundle on $X$.
For a Mukai vector $v$,
$M_L(v)$ is the moduli space of stable sheaves $E$ of 
$v(E)=v$ with respect to $L$. 
For a primitive Mukai vector $v$,
$M_L(v)$ is smooth and projective,
if $L$ is general (cf. [Y2]).

If $\NS(X)={\Bbb Z}$, then
for a coherent sheaf $E$ on $X$,
we set
\begin{equation}
\deg(E):=\frac{(c_1(E),c_1(L))}{(c_1(L)^2)} \in {\Bbb Z},
\end{equation}
where $L$ is the ample generator of $\NS(X)$.

\section{A special case}

\subsection{Preliminaries}

We start with original Fourier-Mukai functor.
So we assume that $X$ is an abelian surface.
Let $\widehat{X}$ be the dual of $X$ and ${\cal P}$ the Poincar\'{e}
line bundle on $X \times \widehat{X}$.
We denote the projections
$X \times \widehat{X} \to X$
(resp. $X \times \widehat{X} \to \widehat{X}$ )
by $p_X$ (resp. $p_{\widehat{X}}$ ).
Let ${\cal F}:{\mathbf D}(X) \to {\mathbf D}(\widehat{X})$ be the
Fourier-Mukai functor defined by 
\begin{equation}
{\cal F}(x):={\mathbf R}p_{\widehat{X}*}({\cal P} \otimes 
p_X^*(x)), x \in {\mathbf D}(X). 
\end{equation}
Let 
$\widehat{{\cal F}}:{\mathbf D}(\widehat{X}) \to 
{\mathbf D}({X})$ be the
inverse of ${\cal F}$:
\begin{equation}
\widehat{{\cal F}}(y):={\mathbf R}p_{{X}*}({\cal P}^{\vee} \otimes 
p_{\widehat{X}}^*(y))[2], y \in {\mathbf D}(\widehat{X}). 
\end{equation}
By [Y5, Lem. 4.2],
we can define an isometry 
${\cal F}_H:H^{ev}(X,{\Bbb Z}) \to H^{ev}(\widehat{X},{\Bbb Z})$ 
of Mukai lattices by
\begin{equation}
{\cal F}_H(x):=p_{\widehat{X}*}((\ch{\cal P})p_X^*\sqrt{\td_X}
 p_{\widehat{X}}^*\sqrt{\td_{\widehat{X}}}
p_X^*(x)), x \in H^{ev}(X,{\Bbb Z}).
\end{equation}
Then the inverse  
$\widehat{{\cal F}}_H:H^{ev}(\widehat{X},{\Bbb Z}) \to 
H^{ev}({X},{\Bbb Z})$ of ${\cal F}_H$ is given by 
\begin{equation}
\widehat{{\cal F}}_H(y):=p_{{X}*}((\ch{\cal P})^{\vee}
p_X^*\sqrt{\td_X}
 p_{\widehat{X}}^*\sqrt{\td_{\widehat{X}}}
p_{\widehat{X}}^*(y)), y \in H^{ev}(\widehat{X},{\Bbb Z}). 
\end{equation}
By Grothendieck Riemann-Roch theorem,
the following diagram is commutative.
\begin{equation}
\begin{CD}
{\mathbf D}(X) @>{{\cal F}}>> {\mathbf D}(\widehat{X})\\
@V{\sqrt{\td_{X}}\ch}VV @VV{\sqrt{\td_{\widehat{X}}}\ch}V\\
H^{ev}(X,{\Bbb Z}) @>{{\cal F}_H}>> H^{ev}(\widehat{X},{\Bbb Z})
\end{CD}
\end{equation}

For a coherent sheaf $E$ on $X$ (resp. a
coherent sheaf $F$ on $\widehat{X}$),
we set
\begin{equation}
\begin{split}
{\cal F}^i(E):&=H^i({\cal F}(E)),\\
\widehat{{\cal F}}^i(F):&=H^i(\widehat{{\cal F}}(F)).
\end{split}
\end{equation}
If $E$ (resp. $F$)
satisfies $\WIT_i$ with respect to ${\cal F}$
(resp. $\widehat{\cal F}$), then 
we denote ${\cal F}^i(E)$ (resp. $\widehat{\cal F}^i(F)$)
by $\widehat{E}$ (resp. $\widehat{F}$).

We are also interested in the composition of ${\cal F}$ and the
``taking-dual'' functor ${\cal D}_{\widehat{X}}:{\mathbf D}(\widehat{X}) \to
{\mathbf D}(\widehat{X})_{op}$
sending $x \in {\mathbf D}(\widehat{X})$ to 
${\mathbf R}{\cal H}om(x,{\cal O}_{\widehat{X}})$,
where ${\mathbf D}(\widehat{X})_{op}$ is the opposite category of 
${\mathbf D}(\widehat{X})$.
By Grothendieck-Serre duality, 
${\cal G}:=({\cal D}_{\widehat{X}} \circ {\cal F})[-2]$
is defined by
\begin{equation}
{\cal G}(x):={\mathbf R}\Hom_{p_{\widehat{X}}}
(p_X^*(x),{\cal P}^{\vee}), x \in {\mathbf D}(X). 
\end{equation}
Let 
$\widehat{{\cal G}}:{\mathbf D}(\widehat{X})_{op} \to 
{\mathbf D}({X})$ be the
inverse of ${\cal G}$:
\begin{equation}
\widehat{{\cal G}}(y):={\mathbf R}\Hom_{p_{{X}}}
(p_{\widehat{X}}^*(y),{\cal P}^{\vee}), y \in {\mathbf D}(\widehat{X}). 
\end{equation}
For a coherent sheaf $E$ on $X$ (resp. a
coherent sheaf $F$ on $\widehat{X}$),
we set
\begin{equation}
\begin{split}
{\cal G}^i(E):&=H^i({\cal G}(E)),\\
\widehat{{\cal G}}^i(F):&=H^i(\widehat{{\cal G}}(F)).
\end{split}
\end{equation}
Then there are spectral sequences
\begin{equation}\label{eq:spectral3}
E_2^{p,q}=\widehat{{\cal G}}^p({\cal G}^{-q}(E)) 
\Rightarrow
\begin{cases}
 E, \;p+q=0\\
0, \text{ otherwise},
\end{cases}
\end{equation}
\begin{equation}\label{eq:spectral4}
E_2^{p,q}={\cal G}^p(\widehat{{\cal G}}^{-q}(F)) 
\Rightarrow
\begin{cases}
 F,\; p+q=0\\
0, \text{ otherwise}.
\end{cases}
\end{equation}
In particular
\begin{equation}\label{eq:vanish2}
\begin{cases}
{\cal G}^p(\widehat{\cal G}^0(F))=0,\;p=1,2,\\
{\cal G}^p(\widehat{\cal G}^2(F))=0,\;p=0,1.
\end{cases}
\end{equation}
If $E$ (resp. $F$)
satisfies $\WIT_i$ with respect to ${\cal G}$
(resp. $\widehat{\cal G}$), then 
we denote ${\cal G}^i(E)$ (resp. $\widehat{\cal G}^i(F)$)
by $\widehat{E}$ (resp. $\widehat{F}$).

\subsection{The case of $\langle v,1 \rangle<0$.}

We assume that $\NS(X)={\Bbb Z} L$,
where $L$ is an ample generator.
Then the dual of $X$ also satisfies the same condition.
We set $\hat{L}:=\det(-{\cal F}(L))$.
Then $\hat{L}$ is the ample generator of $\NS(\widehat{X})$.
For $v=r+dc_1(L)+a \omega \in H^{ev}(X,{\Bbb Z})$,
${\cal F}_H(v)=a-d c_1(\hat{L})+r \omega$.
In this and the next subsections, we shall consider functors
${\cal F}$ and ${\cal G}$.
In this subsection, we treat the case of $\langle v,1 \rangle<0$.
We first treat the functor ${\cal G}$.
\begin{prop}\label{lem:fourier4}
Let $E$ be a $\mu$-stable sheaf of 
Mukai vector $v(E)=r+c_1(L)+a \omega$.
If $a>0$, then $E$ satisfies $\WIT_2$ with respect to ${\cal G}$
and $\widehat{E}$ is 
a $\mu$-stable sheaf of $v(\widehat{E})=a+c_1(\hat{L})+r \omega$.
In particular, ${\cal G}$ induces an isomorphism
$M_L(v) \to M_{\hat{L}}({\cal F}_H(v)^{\vee})$.
\end{prop}
\begin{pf}
(1) $E$ satisfies $\WIT_2$.
By the stability of $E$, ${{\cal G}}^0(E)=0$.
Hence we shall prove that ${{\cal G}}^1(E)=0$.
We first show that $\Ext^1(E, {\cal P}_x^{\vee}) =0$
except for finitely many points $x \in \widehat{X}$.
Suppose that $\Ext^1(E, {\cal P}_x^{\vee})\ne  0$
for distincts points $x_1,x_2,\dots,x_n$.
By [Y4, Lem. 2.1], we get a $\mu$-stable extension sheaf $G$:
\begin{equation}
0 \to \oplus_{i=1}^n {\cal P}_{x_i}^{\vee} \to G
\to E \to 0.
\end{equation}
Since $v(G)=v(E)+n$,
we see that
$\langle v(G)^2 \rangle=\langle v(E)^2 \rangle-2na$.
Hence $n$ must satisfy the inequality
 $n \leq \langle v(E)^2 \rangle/2a$.

By using this, we prove that ${{\cal G}}^1(E)=0$.
By base change theorem,
${{\cal G}}^1(E)$ is of dimension $0$.
Hence we show that $\widehat{{\cal G}}^2({{\cal G}}^1(E))=0$.
Since ${{\cal G}}^0(E)=0$,
$\widehat{{\cal G}}^0({{\cal G}}^0(E))=0$.
By using the spectral sequence \eqref{eq:spectral3}, 
we conclude that $\widehat{{\cal G}}^2({{\cal G}}^1(E))=0$.

(2) 
$\widehat{G}$ is torsion free.\footnote{
This claim also follows from the proof of
base change theorem.} 
Indeed, let $T$ be the torsion submodule of $\widehat{E}$.
Since $\widehat{E}$ is locally free in codimension 1,
$T$ is of dimension 0.
Hence $T$ satisfies $\IT_2$ and
$\deg(\widehat{{\cal G}}^2(T))=0$.
Since $\widehat{{\cal G}}^2(T)$ is a quotient of $E$,
we get a contradiction.

(3)
$\widehat{E}$ is $\mu$-stable.
If $\widehat{E}$ is not $\mu$-stable, then
there is an exact sequence
\begin{equation}
0 \to A \to \widehat{E} \to B \to 0,
\end{equation} 
where $B$ is a $\mu$-stable sheaf of 
$\deg(B) \leq 0$.
Then we get 
\begin{align}
\widehat{{\cal G}}^0(B) &=0,\\
\widehat{{\cal G}}^1(B) &=\widehat{{\cal G}}^0(A),\\
\end{align}
and an exact sequence
\begin{equation}
0 \to \widehat{{\cal G}}^1(A) \to \widehat{{\cal G}}^2(B) \to
E \to \widehat{{\cal G}}^2(A) \to 0.
\end{equation}
If $B \ne {\cal P}_x^{\vee}$ for any $x \in X$,
then $\widehat{{\cal G}}^2(B)=0$.
Hence $B$ satisfies $\WIT_1$.
By \eqref{eq:vanish2},
${\cal G}^1(\widehat{{\cal G}}^2(B))=
{\cal G}^1(\widehat{{\cal G}}^0(A))=0$.
Hence $B=0$, which is a contradiction.
If $B ={\cal P}_x^{\vee}$ for some $x \in X$,
then $\widehat{{\cal G}}^1(B)=0$ and
$\widehat{{\cal G}}^2(B) \cong {\Bbb C}_x$.
Hence $\widehat{{\cal G}}^0(A)=0$ and
$\widehat{{\cal G}}^1(A)=\widehat{{\cal G}}^2(B)= {\Bbb C}_x$.
So we get $ {\cal G}^2(\widehat{{\cal G}}^1(A)) \ne 0$,
which contradicts \eqref{eq:spectral4}.
\end{pf}
By the proof of this proposition, $E$ satisfies $\IT_2$ for ${\cal G}$
if $a>\langle v^2 \rangle/2$.
Hence $E$ also satisfies $\IT_0$ for ${\cal F}$.
Since ${\cal F}^0(E) \cong {\cal G}^2(E)^{\vee}$,
${\cal F}^0(E)$ is also $\mu$-stable.
Thus we get the following.
\begin{cor}\label{lem:fourier1}
Let $E$ be a $\mu$-stable sheaf of 
Mukai vector $v(E)=r+c_1(L)+a \omega$.
We assume that $a>\langle v^2 \rangle/2$.
Then $E$ satisfies $\IT_0$ with respect to ${\cal F}$
and $\widehat{E}$ is 
a $\mu$-stable vector bundle of $v(\widehat{E})=a-c_1(\hat{L})+r \omega$.
\end{cor}
Since $\widehat{{\cal F}}=\widehat{{\cal G}} \circ {\cal D}_{\widehat{X}}$,
we also obtain the following.
\begin{cor}\label{lem:fourier2}
Let $E$ be a $\mu$-stable vector bundle of 
Mukai vector $v(E)=r-c_1(\hat{L})+a \omega$ on $\widehat{X}$.
We assume that $a>0$.
Then $E$ satisfies $\WIT_2$ with respect $\widehat{{\cal F}}$
and $\widehat{E}$ is 
$\mu$-stable.
\end{cor}
In the same way as in the proof of [Y4, Thm. 3.6],
we get another proof of Dekker's result [D, Thm. 5.8]. 
\begin{thm}\label{thm:another}
Let $X$ be an arbitrary abelian surface.
Let $v=r+\xi+a \omega, \xi \in H^2(X,{\Bbb Z})$ be a Mukai vector such that
$r+\xi$ is primitive.
Then $M_L(v)$ is deformation equivalent to
$\widehat{X} \times \Hilb_X^{\langle v^2 \rangle/2}$
for a general ample divisor $L$.
\end{thm}

\subsection{The case of $\langle v,1 \rangle>0$.}
As in the previous subsection, we assume that $\NS(X)={\Bbb Z}L$.

\begin{prop}\label{lem:fourier3}
Let $E$ be a $\mu$-stable sheaf of 
Mukai vector $v(E)=r+c_1(L)+a \omega$.
We assume that $a<0$.
Then $E$ satisfies $\WIT_1$ and $\widehat{E}$ is 
a $\mu$-stable sheaf of $v(\widehat{E})=-a+c_1(\hat{L})-r \omega$.
In particular,
Fourier-Mukai functor induces an isomorphism
$M_L(v) \to M_{\hat{L}}(-{\cal F}_H(v))$.
\end{prop}
\begin{pf}
(1)
$E$ satisfies $\WIT_1$.
We first show that $H^0(X,E \otimes {\cal P}_x) =0$
except for finitely many points $x \in \widehat{X}$.
Suppose that $k_i:=h^0(X,E \otimes {\cal P}_{x_i}) \ne 0$
for distincts points $x_1,x_2,\dots,x_n$.
We shall consider the evaluation map
\begin{equation}
\phi:\oplus _{i=1}^n {\cal P}_{x_i}^{\vee} \otimes 
H^0(X,E \otimes {\cal P}_{x_i}) \to E.
\end{equation}
We assume that $\sum_i k_i >r$,
that is, $\rk(\oplus _{i=1}^n {\cal P}_{x_i}^{\vee} \otimes 
H^0(X,E \otimes {\cal P}_{x_i}))>\rk(E)$.
By the proof of [Y4, Lem. 2.1],
$\phi$ is surjective in codimension 1
and $\ker \phi$ is $\mu$-stable.
We set $b:=\dim(\coker \phi)$.
Then $v(\ker \phi)=\sum_{i=1}^n k_i-(v(E)-b \omega)$.
Since $\sum_i k_i>r$, we get
\begin{equation}
\begin{split}
\langle v(\ker \phi)^2 \rangle &=
\langle v(E)^2 \rangle +2a \sum_i k_i-2b \sum_i k_i+2br\\
& \leq \langle v(E)^2 \rangle +2a \sum_i k_i.
\end{split}
\end{equation}
Since $\langle v(\ker \phi)^2 \rangle \geq 0$, we get
$\sum_i k_i \leq \langle v(E)^2 \rangle/(-2a)$.
In particular, ${\cal F}^0(E)$ is a torsion sheaf.
Since $E$ is a torsion free sheaf on the integral scheme
$\Supp(E)$, ${\cal F}^0(E)$ is torsion free.
Hence we get ${\cal F}^0(E)=0$.
By the stability of $E$, ${\cal F}^2(E)=0$.
Therefore $E$ satisfies $\WIT_1$.

(2)
$\widehat{E}$ is torsion free.
Let $T$ be the torsion subsheaf of $\hat{E}$.
By (1), $T$ is of dimension 0.
Since $\hat{E}$ satisties $\WIT_1$ and $T$ satisfies
$\IT_0$, $T$ must be 0.

(3)
$\widehat{E}$ is $\mu$-stable.
Assume that $\widehat{E}$ is not $\mu$-stable.
Let $0 \subset F_1 \subset F_2 \subset \dots \subset F_s= \widehat{E}$
be the Harder-Narasimhan filtration of $\widehat{E}$.
We shall choose the integer $k$ which satisfies 
$\deg(F_i/F_{i-1})>0, i \leq k$ and
$\deg(F_i/F_{i-1}) \leq 0, i > k$.
We shall prove that
${\cal F}^2(F_k)=0$ and ${\cal F}^0(\widehat{E}/F_k)=0$.
Since $\deg(F_i/F_{i-1})>0, i \leq k$,
semi-stability of $F_i/F_{i-1}$ 
implies that ${\cal F}^2(F_i/F_{i-1})=0, i \leq k$.
Hence ${\cal F}^2(F_k)=0$.
On the other hand, we also see that
${\cal F}^0(F_i/F_{i-1}), i>k$ is of dimension 0.
Since $F_i/F_{i-1}$ is torsion free,
${\cal F}^0(F_i/F_{i-1})=0, i>k$.
Hence we conclude that ${\cal F}^0(\widehat{E}/F_k)=0$.

So $F_k$ and $\widehat{E}/F_k$ satisfy $\WIT_1$ and
we get an exact sequence 
\begin{equation}
0 \to {\cal F}^1(F_k) \to E \to {\cal F}^1(\widehat{E}/F_k) \to 0.
\end{equation}
Since $\deg({\cal F}^1(F_k))=\deg(F_k)>0$,
$\mu$-stability of $E$ implies that
$\deg({\cal F}^1(F_k))=1$ and $\rk({\cal F}^1(F_k))=\rk(E)$.
Thus ${\cal F}^1(\widehat{E}/F_k)$ is of dimension 0.
Then ${\cal F}^1(\widehat{E}/F_k)$ satisfies $\IT_0$,
which is a contradiction.
\end{pf}

\section{More general cases}
In this section, we treat more general cases.
Let $(X,L)$ be a polarized abelian (or K3) surface of
$(L^2)=2r_0k$, where $r_0$ and $k$ are positive integers of
$(r_0,k)=1$.
We assume that $\NS(X)={\Bbb Z}L$.\footnote{Under this assumption,
every simple vector bundle of isotropic primitive Mukai vector
is stable (cf. [Mu3]).}
We set $v_0:=r_0+d_0 c_1(L)+d_0^2k \omega$,
where $d_0$ is an integer of $(r_0,d_0)=1$.
Then $\langle v_0^2 \rangle =0$.
So $Y:=M_L(v_0)$ is an abelian (or K3) surface.
Since $X$ and $Y$ are isogenous, $\NS(Y) \cong {\Bbb Z}$. 
Since $(r_0,d_0^2k)=1$,
there is a universal family ${\cal E}$ on $X \times Y$.
We assume that 
\begin{itemize}
\item
${\cal E}$ is locally free.\footnote{
If ${\cal E}$ is not locally free,
then Fourier-Mukai functor is the same as reflection
functor [Mu3]. This case was treated in [Mr] and [Y4].}
\end{itemize}
We set
\begin{equation}
\hat{L}:=\det(p_{Y!}({\cal E} \otimes {\cal O}_L(kr_0-2kd_0)))^{\vee}.   
\end{equation}
Then $\hat{L}$ is a primitive ample line bundle of 
$(c_1(\hat{L})^2)=(c_1(L)^2)$.
Indeed, $v_0^{\vee}$ and $v({\cal O}_L(kr_0-2kd_0))=c_1(L)-2kd_0 \omega$
generate $(v_0^{\vee})^{\perp}$.
Since 
\begin{equation}
\theta_{v_0^{\vee}}:(v_0^{\vee})^{\perp}/{\Bbb Z}v_0^{\vee}
\to H^2(Y,{\Bbb Z})
\end{equation}
is an isomorphism,
$c_1(\hat{L})$ is primitive.
By Donaldson, $\hat{L}$ is ample.
Hence $\hat{L}$ is a primitive ample line bundle.
Since ${\cal F}^{v_0}({\cal O}_L(kr_0-2kd_0))=
c_1(\hat{L})+b \widehat{\omega}, b \in {\Bbb Z}$
and ${\cal F}^{v_0}_H$ is an isometry of
Mukai lattice,
$(c_1(\hat{L})^2)=
\langle v({\cal F}^{v_0}({\cal O}_L(kr_0-2kd_0)))^2 \rangle=
\langle v({\cal O}_L(kr_0-2kd_0))^2 \rangle
=(c_1(L)^2 )$.

Let ${\cal F}^{v_0}:{\mathbf D}(X) \to {\mathbf D}(Y)$ 
be the Fourier-Mukai functor defined by ${\cal E}$ 
and ${\cal F}_H^{v_0}:H^{ev}(X,{\Bbb Z}) \to H^{ev}(Y,{\Bbb Z})$
the induced isometry.

\begin{lem}\label{lem:image}
Let $d_1$ and $l$ be integers which satisfy
$d_1(kd_0)-lr_0=1$.
Then replacing ${\cal E}$ by ${\cal E}\otimes p_Y^* N$, $N \in \Pic(Y)$,
we get 
\begin{equation}
\begin{cases}
{\cal F}^{v_0}_H(1)=d_0^2 k+d_0 l c_1(\hat{L})+l^2 r_0 \widehat{\omega}\\
{\cal F}^{v_0}_H(c_1(L))=2d_0 k r_0+(2d_0 kd_1-1)c_1(\hat{L})+
(2d_0 k^2 d_1^2-2d_1 k) \widehat{\omega}\\
{\cal F}^{v_0}_H(\omega)=r_0+d_1 c_1(\hat{L})+d_1^2 k \widehat{\omega},
\end{cases}
\end{equation}
where $\widehat{\omega}$ is the fundamental class of $Y$.
\end{lem}
\begin{pf}
We set 
\begin{equation}
\left\{
\begin{split}
&[{\cal F}_H^{v_0}(1)]_1 =a c_1(\hat{L}),\\
&[{{\cal F}_H^{v_0}(c_1(L))}]_1 =b c_1(\hat{L}),\\
&[{\cal F}_H^{v_0}(\omega)]_1 =c c_1(\hat{L}).
\end{split}
\right.
\end{equation}
It is easy to see that $[{\cal F}_H^{v_0}(v_0^{\vee})]_1=0$.
Hence we get the relation
\begin{equation}
r_0a-d_0 b+d_0^2kc=0.
\end{equation}
Since $(r_0,d_0)=1$,
$b \equiv d_0kc \mod r_0$.
By the definition of $\hat{L}$,
$-b+2kd_0c=1$.
Hence we get $kd_0c \equiv 1 \mod r_0$.
So replacing ${\cal E}$ by
${\cal E} \otimes \hat{L}^{\otimes((d_1-c)/r_0)}$,
we may assume that
\begin{equation}
\left\{
\begin{split}
& [{\cal F}_H^{v_0}(c_1(L))]_1=(2kd_0d_1-1) c_1(\hat{L}),\\
& [{\cal F}_H^{v_0}(\omega)]_1=d_1 c_1(\hat{L}).
\end{split}
\right.
\end{equation}
Since ${\cal F}_H^{v_0}$ is an isometry,
$\langle {\cal F}_H^{v_0}(\omega)^2 \rangle=
\langle \omega^2 \rangle=0$.
Hence we get
\begin{equation}
{\cal F}_H^{v_0}(\omega)=r_0+d_1 c_1(\hat{L})+d_1^2k \widehat{\omega}.
\end{equation}
Since ${\cal E}$ is a universal family of stable sheaves of
Mukai vector $v_0$,
we get the following relations:
\begin{equation}
\begin{cases}
{\cal F}_H^{v_0}(v_0^{\vee})=\widehat{\omega}\\
{\cal F}_H^{v_0}(\omega)=r_0+d_1 c_1(\hat{L})+d_1^2k \widehat{\omega}\\
{\cal F}_H^{v_0}(-c_1(L)+2kd_0 \omega)=x+c_1(\hat{L})+y \widehat{\omega},
\end{cases}
\end{equation}
where $x,y \in {\Bbb Z}$.
Since ${\cal F}_H^{v_0}$ is an isometry,
we see that $x=0$ and $y=-2kd_1$.
Hence we get our lemma.
\end{pf}
In order to generalize Proposition \ref{lem:fourier4},
and \ref{lem:fourier3}, let us introduce some notations.
Let $G$ be a locally free sheaf on $X$.
For a torsion free sheaf $E$ on $X$, we define
\begin{equation}
\begin{split}
\rk_G(E):& =\rk(E \otimes G^{\vee}),\\
\deg_G(E):&=\deg (E \otimes G^{\vee}),\\
\mu_G(E):& =\frac{\deg_G(E)}{\rk_G(E)}.
\end{split}
\end{equation}
For $x \in {\mathbf D}(X)$ (resp. $v(x) \in H^{ev}(X,{\Bbb Z})$),
we can also define $\rk_G(x)$ and $\deg_G(x)$
(resp. $\rk_G(v(x))$ and $\deg_G(v(x))$ .

Then we see that 
\begin{equation}
\begin{split}
\mu_G(E) & =\frac{\deg(E) \rk (G)-\deg(G) \rk(E)}{\rk(G)\rk(E)}\\
&=\mu(E)-\mu(G).
\end{split}
\end{equation}  
Hence $E$ is $\mu$-stable if and only if
\begin{equation}
\mu_G(F)<\mu_G(E)
\end{equation}
for any subsheaf $F \subset E$ of $\rk(F)<\rk(E)$.
Assume that $\deg_G(E)=1$. Then it is easy to see that
$E$ is $\mu$-stable if and only if
$\deg_G(F) \leq 0$ for any subsheaf $F \subset E$ of $\rk(F)<\rk(E)$.  
\begin{lem}\label{lem:deg}
We choose points $s \in X$ and $t \in  Y$. 
We set 
\begin{equation}
\begin{split}
G_1 &:={\cal E}_{|X \times \{t \}}^{\vee},\\
G_2 &:={\cal E}_{|\{s\} \times Y}.
\end{split}
\end{equation}
Then for a Mukai vector $v$,
\begin{equation}
\deg_{G_1}(v)=-\deg_{G_2}({\cal F}^{v_0}_H(v))=
\deg_{G_2^{\vee}}({\cal F}^{v_0}_H(v)^{\vee}).
\end{equation}
\end{lem}
\begin{pf}
We set 
\begin{equation}
\begin{split}
v&=r+d c_1(L)+a \omega,\\
{\cal F}^{v_0}_H(v)&=r'+d' c_1(L)+a'\omega.
\end{split}
\end{equation}
It is sufficient to prove that
$r'd_1-d'r_0=dr_0+rd_0$.
This follows from the following relations which comes from
Lemma \ref{lem:image}: 
\begin{equation}
\begin{cases}
r'=r(d_0^2 k)+d(2d_0 r_0 k)+a r_0,\\
d'=r(d_0 l)+d(2d_0 d_1 k-1)+ad_1.
\end{cases}
\end{equation}
\end{pf}
Due to this lemma, we can use the same 
arguments as in Propositions \ref{lem:fourier4}
and \ref{lem:fourier3}.
Hence we get the following theorem.
\begin{thm}\label{thm:fourier}
Keep the notations as above.
Let $v:=r+d c_1(L)+a \omega$ be a Mukai vector
of $d r_0+r d_0=1$.
\begin{enumerate}
\item
If $-\langle v,v_0^{\vee} \rangle>0$, then
the composition of Fourier-Mukai functor and
``taking-dual'' functor ${\cal D}_Y$ induces an isomorphism
$M_L(v) \to M_{\hat{L}}({\cal F}^{v_0}_H(v)^{\vee})$.
\item
If $\langle v,v_0^{\vee} \rangle>0$, then Fourier-Mukai functor induces
an isomorphism
$M_L(v) \to M_{\hat{L}}(-{\cal F}^{v_0}_H(v))$.
\end{enumerate}
\end{thm}
\begin{ex}
Keep the notations as above.
We set 
\begin{equation}
\begin{cases}
d_0=-(r_0-1),\\
r=d=1,\\
k=-n+s r_0,\\
a=(r_0^2-1)s-r_0n,
\end{cases}
\end{equation}
where $s>0$, $s r_0>n > 0$ and $(r_0,n)=1$.
Then $\langle v^2 \rangle=2s$ and
$\langle v,v_0^{\vee} \rangle=n>0$.
Applying Theorem \ref{thm:fourier},
we get an isomorphism
$M_L(v) \to M_{\hat{L}}(-{\cal F}_H^{v_0}(v))$.
In particular,
we get another proof of [Y4, Thm. 0.2] and
Theorem \ref{thm:another}.   
\end{ex}

\begin{ex}
We assume that $X$ is a K3 surface and $(L^2)=12$.
We set $v_0:=2-L+3 \omega$ and
$Y:=M_L(v_0)$.
Then $Y$ is a K3 surface of 
$H^2(Y,{\Bbb Z}) \cong v_0^{\perp}/{\Bbb Z}v_0$.
In general, $Y \ne X$.
By Fourier-Mukai functor defined by $v_0$,
we get an isomorphism
$M_L(1+L+3 \omega) \cong M_{\hat{L}}(3-\hat{L}+\widehat{\omega})$.
By reflection $R_{v({\cal O}_Y)}$ ([Mr], [Y4]),
we obtain a birational map 
$M_{\hat{L}}(3-\hat{L}+\widehat{\omega}) \cdots \to
M_{\hat{L}}(1+\hat{L}+3\widehat{\omega})$.
Hence we get a birational map
$\Hilb_X^4 \cdots \to \Hilb_Y^4$.
Computing ample cones, we see that this birational map is the
elementary transformation along ${\Bbb P}^2$-bundle over
$X \times Y$.
In this case, by Torelli Theorem for K3 surfaces,
we see that $\Hilb_X^4 \not \cong \Hilb_Y^4$.
\end{ex}

\section{appendix}

\begin{lem}
Keep the notation in section 2.
Let $E$ be a $\mu$-stable locally free sheaf
of $\deg_{G_1}(E)=0$ and $E \not \in M_L(v_0^{\perp})$.
Then $E$ satisfies $\IT_1$ and
$\widehat{E}$ is a $\mu$-stable locally free sheaf.
\end{lem}
\begin{pf}
Since $E \not \in M_L(v_0^{\perp})$ and $E$ is $\mu$-stable,
we see that $E$ satisfies $\IT_1$.
In the same way as in the proof of Lemma \ref{lem:fourier3},
we see that $\widehat{E}$ is $\mu$-semi-stable.
Assume that there is an exact sequence
\begin{equation}
0 \to F_1 \to \widehat{E} \to F_2 \to 0,
\end{equation}
where $F_2$ is a $\mu$-stable sheaf of 
$\deg_{G_2}(F_2)=0$ and $\rk(F_2)<\rk(\widehat{E})$.
Then $F_2$ satisfies $\IT_1$ and we get an exact
sequence
\begin{equation}
0 \to {\cal F}^1(F_1) \to E \to {\cal F}^1(F_2) \to 
{\cal F}^2(F_1) \to 0.
\end{equation} 
Since ${\cal F}^2(F_1)$ is of dimension 0,
$\deg_{G_1}({\cal F}^1(F_1))=0$.
By the $\mu$-stability of $E$, 
we get ${\cal F}^1(F_1)=0$ or $\rk({\cal F}^1(F_1))=\rk(E)$.
We first assume that ${\cal F}^1(F_1)=0$.
 Since $E$ and ${\cal F}^1(F_2)$ are locally free,
${\cal F}^2(F_1)=0$. Hence $F_1=0$, which is a contradiction.
We next assume that $\rk({\cal F}^1(F_1))=\rk(E)$.
Since ${\cal F}^1(F_2)$ is locally free,
${\cal F}^1(F_1) \to E$ is an isomorphism.
Hence ${\cal F}^1(F_2) \cong {\cal F}^2(F_1)$.
Since ${\cal F}^2(F_1)$ is of dimension 0,
we obtain ${\cal F}^1(F_2)=0$.
Therefore $F_2=0$, which is a contradiction.
\end{pf}

\vspace{1pc}

{\it Acknowledgement.}
I learned Dekker's thesis from
Professor T. Katsura and Professor G. van der Geer,
which motivated me very much.
I would like to thank them very much.
I would also like to thank Max Planck Institut f\"{u}r Mathematik
for support and hospitality.

\end{document}